\documentclass[letterpaper, preprint, paper,11pt]{AAS}	

\usepackage{bm}
\usepackage{amsmath}
\usepackage{amssymb}
\usepackage{subfig}
\usepackage[colorlinks=true, pdfstartview=FitV, linkcolor=black, citecolor= black, urlcolor= black]{hyperref}
\usepackage{overcite}
\usepackage{footnpag}			      	
\usepackage{comment}

\newcommand{\hide}[1]{}  

\newcommand{\STMa}[2]{\bm{\Phi}^{\textbf{#1}}_{\textbf{#2}}}

\newcommand{\STMb}[2]{\bm{\Phi}^{\textbf{#1}}_{\bm{#2}}}

\newcommand{\STMc}[2]{\bm{\Phi}^{\bm{#1}}_{\textbf{#2}}}

\newcommand{\STMd}[2]{\bm{\Phi}^{\bm{#1}}_{\bm{#2}}}

\PaperNumber{25-153}

\begin{document}


\title{Mass-Optimal Low-Thrust Forced Periodic Trajectories in the Earth-Moon CR3BP}


\author{Colby C. Merrill\thanks{Ph.D. Candidate, Sibley School of Mechanical and Aerospace Engineering, Cornell University, Ithaca, NY 14853, USA}, 
Jackson Kulik\thanks{Assistant Professor, Department of Mechanical and Aerospace Engineering, Utah State University, Logan, UT 84322, USA}, 
Matthew J. Bryan\thanks{M.Eng., Sibley School of Mechanical and Aerospace Engineering, Cornell University, Ithaca, NY 14853, USA}, 
Dmitry Savransky\thanks{Associate Professor, Sibley School of Mechanical and Aerospace Engineering, Cornell University, Ithaca, NY 14853, USA}}

\maketitle{}


\begin{abstract} 
In Cislunar space, spacecraft are able to exploit naturally periodic orbits, which provide operational reliability. However, these periodic orbits only exist in a limited volume. Enabled by low-thrust propulsion, spacecraft can produce a greater number of periodic trajectories in Cislunar space. We describe a methodology for producing mass-optimal trajectories that enforce periodic structure in the circular-restricted three body problem and study the thrust-limited reachable set around a reference trajectory. In this study, we find that the thrust-limited mass-optimal reachable set is a superset of the energy-limited energy-optimal reachable set in the $xy$-plane.
\end{abstract}


\section{Introduction}

Spacecraft operating in the Earth-Moon system often orbit on (or in the vicinity of) natural (unforced), periodic or quasi-periodic orbits. Periodic orbits are  naturally bounded and so it is operationally convenient and straightforward to perform stationkeeping for spacecraft on these orbits. However, these orbits only exist in specific areas in Cislunar space and so they do not offer significant coverage of the Cislunar volume. Thus, these regions of natural periodic orbits may limit the operational possibilities for spacecraft and limit the ability for a spacecraft to perform certain tasks. Many modern spacecraft control their trajectory with onboard, low-thrust propulsion systems. Spacecraft with low-thrust propulsion capabilities may exploit trajectories that expand the region of operational capability of Cislunar space and simultaneously achieve the attractive operational attributes of periodic orbits. In this work, we refer to these orbits as ``forced periodic trajectories.'' 

To analyze trajectories in Cislunar space, a common treatment is to reduce the dynamical complexity of the system to a model known as the circular restricted three body problem (CR3BP). This model includes the gravitational effects of the Earth and Moon and assumes that both bodies are on mutually circular orbits. In this model, the locations of five equilibirum points (i.e., Lagrange points) may be solved for along with a number of periodic orbit families. A perturbation applied to the CR3BP (e.g., a constant acceleration provided by low thrust propulsion), will shift the five equilibrium points \cite{McInnes1994, Morimoto2007, Baig2008, Cox2019, Cox2020}. Because the equilibrium points have been shifted, the periodic structures around the equilibria are shifted as well \cite{Morimoto2006}. Previous studies have primarily focused on constant, low thrust trajectories \cite{Cox2019, Cox2020, Morimoto2006, DeLeo2022, Baig2009} or optimal control in proximity of the shifted equilibria \cite{Tsuruta2024}. Other related work exists that focuses on forced circumnavigation or controlled loitering trajectories relative to a spacecraft on some reference orbit in the CR3BP \cite{Sandel2024}. In previous work \cite{Merrill2024}, we studied the reachable set of forced periodic structures under energy-optimality in proximity to a reference orbit. In this work, we focus on mass-optimal forced periodic structures under realistic operational constraints for a spacecraft in Cislunar space. These orbits can also be viewed as forced periodic orbits relative to a chief satellite on a naturally periodic orbit \cite{Sandel2024} but we specifically discuss them in terms of Cislunar trajectories. 

In previous work \cite{Merrill2024}, we applied a linear analysis to find a closed-form solution to the boundary value problem and studied the set of forced periodic trajectories around a reference trajectory that satisfy an energy cost constraint. Although this approach lends itself to studying a set of energy-limited forced periodic trajectories, it also has a number of shortcomings that we address in this paper. First, an energy-optimal solution is not necessarily ideal in spacecraft operations. Instead, a mass-optimal thrust-limited solution is preferred, as these solutions coincide with spacecraft operational objectives and constraints. Second, the energy-limiting approach does not enforce thrust limits in the analysis, which results in expansion of the extents of the reachable set in certain directions. To enforce constraints and relate forced periodic trajectories to operational conditions for a spacecraft, we study the thrust-limited reachable set around a reference trajectory. 


\section{Methods}

\subsection{Dynamics}
We assume a dynamical system of the form
\begin{equation}
    \frac{\mathrm{d}\mathbf{x}}{\mathrm{d}t}=\mathbf{F}(\mathbf{x})+\begin{bmatrix}\mathbf{0}_{3\times1}\\
        \mathbf{u}
    \end{bmatrix}
    \label{eq:dynamics}
\end{equation}
where the state vector $\mathbf{x}\in\mathbb{R}^6$ is defined by stacking the position and velocity vectors $\mathbf{x}=[\mathbf{r}^T, \mathbf{v}^T]^T$ and $\mathbf{u}$ is the control acceleration vector. In our system, $\mathbf{F}(\mathbf{x})$ gives the natural dynamics for the third body in the canonical rotating frame of the CR3BP
\begin{equation}
    \mathbf{F}(\mathbf{x}) =
    \begin{bmatrix}
        v_{x}\\ v_{y}\\ v_{z}\\ 
        2v_{y} + x - (1-\mu^*)\frac{x+\mu^*}{R_{3/1}^3} - \mu^*\frac{x-1+\mu^*}{R_{3/2}^3}\\
        -2v_{x} + y - (1-\mu^*)\frac{y}{R_{3/1}^3} - \mu^*\frac{y}{R_{3/2}^3}\\
        -(1-\mu^*)\frac{z}{R_{3/1}^3} - \mu^*\frac{z}{R_{3/2}^3}
    \end{bmatrix}
\end{equation}
where $x$, $y$, and $z$ are the components of the spacecraft's position vector, $\mathbf{r}$, and $v_x$, $v_y$, and $v_z$ are the components of the spacecraft's velocity vector, $\mathbf{v}$. In this frame, the $\hat{\mathbf{r}}_x$ direction points from the barycenter of the primary and secondary bodies toward the secondary body, the $\hat{\mathbf{r}}_z$ direction is parallel to the direction of the angular momentum vector of the two larger masses, and the $\hat{\mathbf{r}}_y$ direction completes the right-hand coordinate system. The distances of the third body with respect to the primary and secondary are defined as $R_{3/1}$ and $R_{3/2}$, respectively, and are evaluated as
\begin{align}
    R_{3/1} &= \sqrt{(x+\mu^*)^2 + y^2 + z^2}\\
    R_{3/2} &= \sqrt{(x-1+\mu^*)^2 + y^2 + z^2}
\end{align}
where $\mu^* = m_2/(m_1+m_2)$ is the mass parameter of the system. Note that to improve numerical precision in our computation, components in this frame and masses are measured in canonical units [TU], [DU], and [MU] which satisfy the relationships
\begin{align}
    1 \: \text{[TU]} &\equiv \sqrt{\frac{R_{2/1}^3}{G(m_1+m_2)}}\\
    1 \: \text{[DU]} &\equiv R_{2/1}\\
    1 \: \text{[MU]} &\equiv m_1+m_2
\end{align}
where $R_{2/1}$ is the distance between the primary and secondary bodies. This allows us to express the gravitational constant G as
\begin{equation}
    G = 1 \left[\frac{\text{DU}^3}{\text{TU}^2\text{MU}}\right]
\end{equation}

\subsection{Optimization}
The optimal control from one state to another is given by solving a two-point boundary value problem associated with a system of ordinary differential equations. In the indirect approach, these equations have twice as many dimensions as the state of the original system and are given by
\begin{align}
    \frac{\mathrm{d}\mathbf{x}}{\mathrm{d}t}&=\mathbf{F}(\mathbf{x})+\begin{bmatrix}
        \mathbf{0}_{3\times1}\\\mathbf{u}\end{bmatrix}\\
    \frac{\mathrm{d}\bm{\lambda}}{\mathrm{d}t}&=-\left(\frac{\partial \mathbf{F(x)}}{\partial \mathbf{x}}\right)^T \bm{\lambda}    \\
    \mathbf{u} &= -\bm{\lambda_v}
\end{align}
where $\bm{\lambda_v}$ is the velocity costate vector given by the last three elements of the costate vector \cite{Bryson2018}. In the direct approach, the only differential equations are those listed in Equation \ref{eq:dynamics}. The additional constraints that we enforce are
\begin{align}
    \mathbf{x}(t_0) &= \mathbf{x}_0 \\
    \mathbf{x}(t_f) &= \mathbf{x}(t_0) \\
    ||\mathbf{u}(t)|| &\le u_{max}
\end{align}
where $\mathbf{x}_0$ is some selected initial state for the spacecraft and $u_{max}$ is the maximum thrust output from the spacecraft. To find energy-optimal trajectories, we minimize a cost function of the form
\begin{equation}
    J_E = \int_{t_0}^{t_f} \frac{1}{2}||\mathbf{u}(t)||^2 \mathrm{d}t
    \label{eq:J_energy}
\end{equation}
such that 
\begin{align}
\arg \min_{\mathbf{u}(t)} J_E = \mathbf{u}^*_E
\end{align}
where $\mathbf{u}^*_E$ is the control history of an energy-optimal trajectory. Similarly, to find mass-optimal trajectories, we minimize the cost function 
\begin{equation}
    J_M = \int_{t_0}^{t_f} ||\mathbf{u}(t)|| \mathrm{d}t
    \label{eq:J_mass}
\end{equation}
such that 
\begin{align}
\arg \min_{\mathbf{u}(t)} J_M = \mathbf{u}^*_M
\end{align}
where $\mathbf{u}^*_M$ is the control history of a mass-optimal trajectory. 

Here, we optimize the full nonlinear problem with the use of the Astrodynamics Software and Science Enabling Toolkit (ASSET) \cite{Pezent2022}. With ASSET, we can determine the spacecraft's optimal thrust profile given a set of applied constraints and dynamics to the optimal control problem. To handle the optimization, parallel sparse interior-point optimizer (PSIOPT) is used to solve the non-linear programming problem (NLP) with a primal-dual-interior-point method. Important to the form of the optimal control problem we investigate, PSIOPT can parallelize functions for rapid computation. We use a high-order Legendre Gauss Lobatto collocation method for transcription, which has been written in ASSET. Note that this transcription method does not directly control integration error, so ASSET also has a mesh refinement that continuously iterates the transcription to satisfy given error tolerances \cite{Pezent2023}. For a thorough discussion of ASSET and PSIOPT, we refer the reader to chapter 2 of Pezent 2024 \cite{Pezent2024}.

\subsection{Generation of Mass-Optimal Trajectories}
We follow a specific procedure to generate our mass-optimal solutions in this paper. In contrast to previous studies that use smoothing functions or homotopy to achieve mass-optimal results, we are able to find mass-optimal solutions directly from an energy-optimal initial guess without intermediate steps. Our procedure is as follows:
\begin{enumerate}
    \item Select initial conditions for a naturally periodic trajectory in some given three-body system to serve as a reference. If needed, propagate this naturally periodic trajectory and iterate on the initial conditions with differential correction such that orbit closes on itself in a single period.
    \item Use the propagated reference trajectory as the initial guess to solve a thrust-constrained energy-optimal problem with a cost function given by Equation \ref{eq:J_energy}.
    \item Use the output from the energy-optimal problem as the initial guess to solve a constrained mass-optimal problem with a cost function given by Equation \ref{eq:J_mass}. 
    \item Apply mesh refinement to reduce the error in the mass-optimal solution.
    \item Reintegrate the mass-optimal trajectory with the thrust history to verify that the solution satisfies the problem's constraints. 
\end{enumerate}
From our experience, we find that using a third-order LGL collocation method with $\sim$ 50 knot points is able to produce the energy-optimal result. To find the mass-optimal result, we again use a third-order LGL collocation method with $\sim$ 100 knot points. For the mesh refinement step, we set our mesh tolerance to 1E-10. For reintegration, we use an adaptive 8(7) Dormand-Prince integrator with a relative tolerance of 1E-14 and absolute tolerance of 1E-16. 

The reference trajectory used throughout this paper is the same used in Merrill et al.~2024\cite{Merrill2024} and is given by the initial conditions
\begin{equation}
\mathbf{x}_{ref} = 
    \begin{bmatrix}
    1.06315768  \:\: \text{DU}\\  0.000326952322 \:\: \text{DU}\\ -0.200259761 \:\: \text{DU}\\  0.000361619362 \:\: \text{DU/TU}\\ -0.176727245 \:\: \text{DU/TU}\\ -0.000739327422 \:\: \text{DU/TU}
    \end{bmatrix} 
\end{equation}
with a period of 2.085034838884136 TU and a mass constant of 0.01215059. In this study, we assume a spacecraft with an initial mass of 1000 kg and maximum thrust of 50 mN, translating to $u_{max} =5 \times 10^{-5}$ m/s$^2 \approx 0.0184$ DU/TU$^2$.

\subsection{Reachable Set Definitions}

The energy-limited reachable set is defined by 
\begin{equation}
    \left\{ \delta\mathbf{x}_0 \quad \mathrm{s.t.}\quad J_E \leq \frac{1}{2}u_{max}^2(t_f-t_0) \right\}
\end{equation}
where $\delta\mathbf{x}_0$ is the deviation from some reference trajectory's initial state (i.e., $\mathbf{x}_0 = \mathbf{x}_{ref} + \delta\mathbf{x}_0$). This reachable set takes the form of a six-dimensional hyper-ellipsoid \cite{Merrill2024}. The thrust-limited reachable set is defined by 
\begin{equation}
    \left\{ \delta\mathbf{x}_0 \quad \mathrm{s.t.}\quad ||\mathbf{u}^*_M(t)|| \leq u_{max} \quad \mathrm{for}\quad t_0\leq t\leq t_f \right\}
\end{equation}
and cannot be described analytically in the same way as the energy-limited reachable set.  The difference between the sets is central to this analysis. In Kulik et al.~2024\cite{Kulik2024}, the differences between the thrust-limited and energy-limited reachable sets were studied for near-circular orbits. In this analysis, we sample from the boundaries of the reachable set. These definitions have some external constraints as well including the requirement that the forced periodic trajectory has the same period as the reference trajectory. Both the linearized method used in Merrill et al.~2024\cite{Merrill2024} and the method used in this paper allow for this requirement to be relaxed, but that is not the focus of the current work.

\subsection{Energy-Limited Reachable Set}

To obtain an initial analytical understanding of the problem, we first find the reachable set of energy-optimal forced periodic trajectories limited by the total energy cost. We refer the reader to the work of Merrill et al.~2024\cite{Merrill2024} for a full explanation of this process, as we will only summarize it here. 

First, we define the augmented state as
\begin{equation}
\mathbf{y} = \begin{bmatrix}
        \mathbf{x} \\ 
        \bm{\lambda}
    \end{bmatrix} =
    \begin{bmatrix}
        \mathbf{r} \\
        \mathbf{v} \\
        \bm{\lambda}_r \\
        \bm{\lambda_v}
    \end{bmatrix}
\end{equation}
and note that $J_E$ can then be written in terms of the velocity costate vector as
\begin{equation}
    J_E = \frac{1}{2}\int_{t_0}^{t_f} \bm{\lambda_v}^T \bm{\lambda_v} \mathrm{d}t
\end{equation}
The state transition matrix (STM) associated with the augmented state vector and its dynamics yields a linear approximation of perturbations to the final augmented state $\delta\mathbf{y}(t_f)$ at some final time, $t_f$, as a function of deviations in the initial augmented state:
\begin{equation}
\delta\mathbf{y}(t_f) = \begin{bmatrix}
        \delta\mathbf{r}(t_f) \\ \delta\mathbf{v}(t_f) \\ \delta\bm{\lambda}_r(t_f) \\ \delta\bm{\lambda_v}(t_f)
    \end{bmatrix} \approx \STMd{}{}(t_f,t_0) \delta\mathbf{y}_0 = \begin{bmatrix}
      \STMa{x}{x} & \STMb{x}{\lambda}  \\ \STMc{\lambda}{x} & \STMd{\lambda}{\lambda}
    \end{bmatrix} \delta\mathbf{y}_0
\end{equation}
where $\delta\mathbf{y}_0 = \delta\mathbf{y}(t_0)$ is the perturbed initial state and $\bm{\Phi}(t_f,t_0)$ is a time-varying STM associated with the augmented state, reference trajectory, and initial and final times. We adopt the notation
\begin{equation}
    \STMa{b}{a}(t,t_0) = \frac{\partial \mathbf{b}(t)}{\partial \mathbf{a}(t_0)}
\end{equation}
When the time dependence of an STM is omitted in this paper, it indicates that the STM corresponds to a full period (i.e., $\STMd{}{}(t_f,t_0) = \STMd{}{}$). $\delta\bm{\lambda_v}(t)$ at all times may be evaluated as 
\begin{equation}
\delta\bm{\lambda_v}(t) = \STMc{\lambda_v}{y}(t,t_0) \delta\mathbf{y}_0 = 
    \begin{bmatrix}
      \STMc{\lambda_v}{r}(t,t_0) & \STMc{\lambda_v}{v}(t,t_0) & \STMd{\lambda_v}{\lambda_r}(t,t_0) & \STMd{\lambda_v}{\lambda_v}(t,t_0)
    \end{bmatrix} \delta\mathbf{y}_0
\end{equation}
and substituted into the cost function as
\begin{equation}
    J_E = \frac{1}{2}\delta\mathbf{y}_0^T \left(\int_{t_0}^{t_f} \left(\STMc{\lambda_v}{y}(t,t_0)\right)^T \left(\STMc{\lambda_v}{y}(t,t_0)\right) dt\right) \delta\mathbf{y}_0
    \label{eq:dy_cost_integral}
\end{equation}
This linearized analysis allows us to explicitly solve for the initial costates that satisfy the linearized boundary value problem constraints
\begin{equation}
    \delta\bm{\lambda}_0 = \begin{bmatrix}
    -\left(\STMb{x}{\lambda}\right)^{-1} \STMa{x}{x} & \left(\STMb{x}{\lambda}\right)^{-1}  \end{bmatrix} \begin{bmatrix}
    \delta\mathbf{x}_0 \\ \delta\mathbf{x}_f
    \end{bmatrix}
\end{equation}
which then can be used to find the full, initial, augmented state in terms of the boundary conditions
\begin{equation}
    \delta\mathbf{y}_0 = \begin{bmatrix}
    \mathbf{I}_6 & \mathbf{0}_6 \\ -\left(\STMb{x}{\lambda}\right)^{-1} \STMa{x}{x} & \left(\STMb{x}{\lambda}\right)^{-1} \end{bmatrix} \begin{bmatrix}
    \delta\mathbf{x}_0 \\ \delta\mathbf{x}_f
    \end{bmatrix}
\end{equation}
We define the matrix
\begin{equation}
    \mathbf{E} =  \begin{bmatrix} 
    \mathbf{I}_6 & \mathbf{0}_6 \\ -\left(\STMb{x}{\lambda}\right)^{-1} \STMa{x}{x} & \left(\STMb{x}{\lambda}\right)^{-1} 
    \end{bmatrix}^T \int_{t_0}^{t_f} \left(\STMc{\lambda_v}{y}(t,t_0)\right)^T \left(\STMc{\lambda_v}{y}(t,t_0)\right) dt 
    \begin{bmatrix} \mathbf{I}_6 & \mathbf{0}_6 \\ -\left(\STMb{x}{\lambda}\right)^{-1} \STMa{x}{x} & \left(\STMb{x}{\lambda}\right)^{-1} 
    \end{bmatrix}
\end{equation}
which can be used to determine the energy-constrained reachable set for our system \cite{Lee2018, Sun2020}. Substituting in to the cost function, we now have 
\begin{equation}
    J_E = \frac{1}{2}\begin{bmatrix}
    \delta\mathbf{x}_0 \\ \delta\mathbf{x}_f
    \end{bmatrix}^T \mathbf{E}
    \begin{bmatrix}
    \delta\mathbf{x}_0 \\ \delta\mathbf{x}_f
    \end{bmatrix}
\end{equation}
Forced periodic trajectories in the vicinity of a periodic reference orbit satisfy the condition $\delta\mathbf{x}_f=\delta\mathbf{x}_0$, so that the final state is equivalent to the initial state after one period of the reference periodic orbit. To study the set of forced periodic trajectories that require less than some energy limit, we may study the following matrix 
\begin{equation}
    \mathbf{E}^* = \begin{bmatrix}
        \mathbf{I}_6 & \mathbf{I}_6
    \end{bmatrix} \mathbf{E} \begin{bmatrix}
        \mathbf{I}_6 & \mathbf{I}_6
    \end{bmatrix}^T
\end{equation}
so that the linearized cost function for the state-return trajectory beginning and ending at $\delta\mathbf{x}_0$ can then be written as
\begin{equation}
    J_E = \frac{1}{2}\delta\mathbf{x}_0^T \mathbf{E}^* \delta\mathbf{x}_0
\end{equation}
Now assume that ($\gamma_i, \mathbf{w}_i$) is an eigenpair of the matrix $\mathbf{E}^*$ where $\gamma_i$ is some eigenvalue and $\mathbf{w}_i$ is its corresponding eigenvector. Since $\mathbf{E}^*$ is a symmetric positive semi-definite matrix, the set of possible relative states $\delta\mathbf{x}_0$ that cost less than some energy limit $J^*$ to begin and end at under linearized optimal control is given by the hyperellipsoid described by the set
\begin{equation}
    \left\{ \delta\mathbf{x}_0 \quad \mathrm{s.t.}\quad \frac{1}{2}\delta\mathbf{x}_0^T \mathbf{E}^* \delta\mathbf{x}_0 \leq \frac{1}{2}u_{max}^2(t_f-t_0) \right\}
\end{equation}
Thus, we have defined the set of relative states which can be returned to in a period of the reference orbit. The ellipsoid described here is in 6-dimensional position and velocity space. 

\subsection{Particle Swarm Optimization}

Finding the reachable set of mass-optimal trajectories poses a challenge to gradient-based optimization procedures that require smoothness and continuity. By searching for solutions at the boundary between feasible and infeasible trajectories, we push the solutions up to the barrier of the convergence space. Furthermore, given the setup of the problem, we cannot readily re-parameterize the problem to explicitly search for the surface of this reachable set while maintaining all constraints and nonlinearity. To address these problems, we currently use an accelerated particle swarm optimization (PSO) algorithm\cite{Kennedy1995, Gandomi2013} to find the maximum deviations in a specified direction. 

\begin{figure}[ht]
    \centering
    \includegraphics[width=1.0\linewidth]{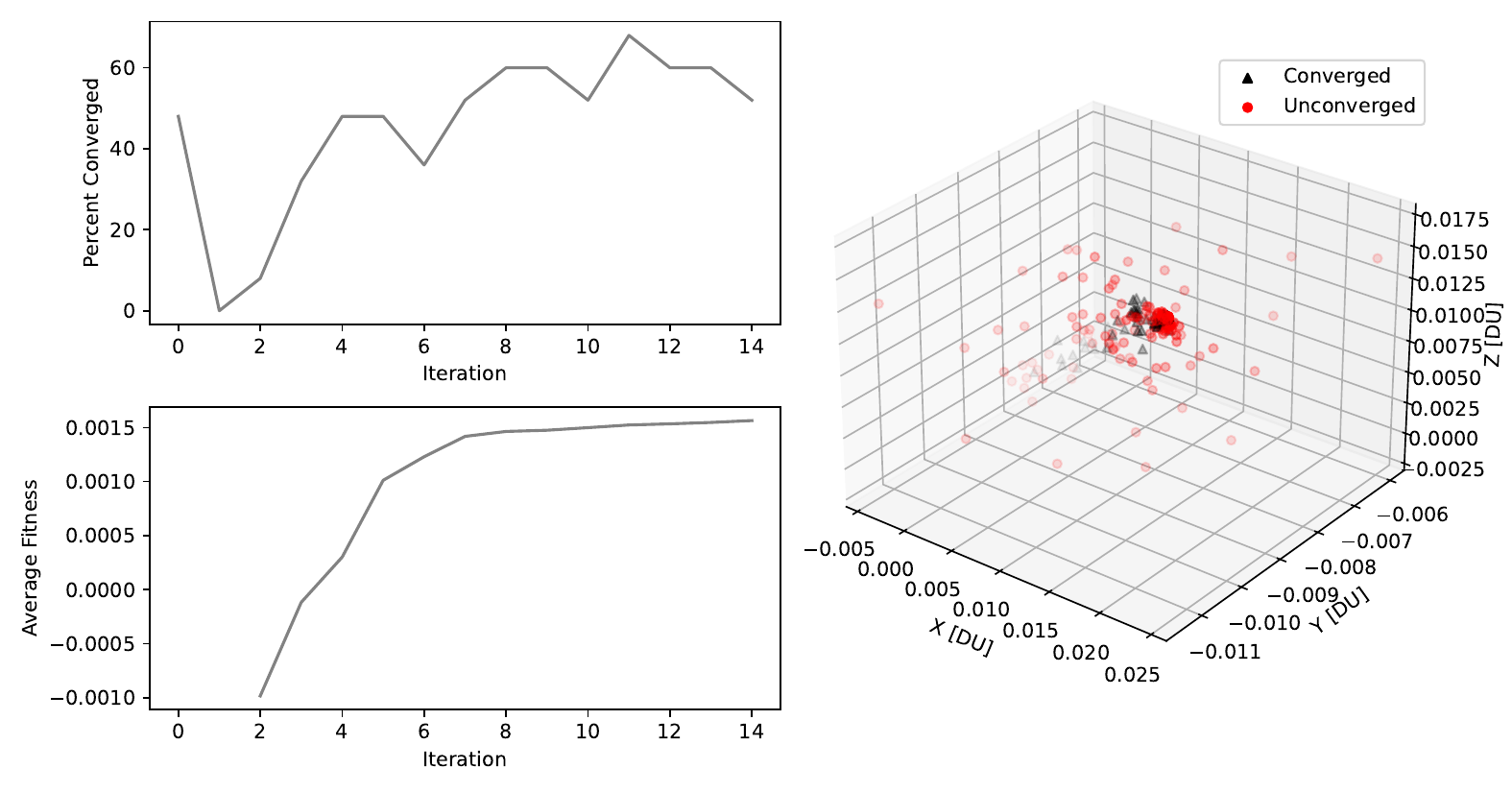}
    \caption{An example of the PSO algorithm.}
    \label{fig:pso_example}
\end{figure}

At the beginning of each PSO run, we initialize the particles and center them around a known solution in the proximity of some base particle. This base particle can be a known solution to a similar optimization problem (i.e., PSO can be used similar to a continuation method) or a vector of zeroes. The components of the state are altered from the base particle by adding in points from a normal distribution. This distribution should have standard deviations on the scale of the search space of interest. For this problem, each component of the state is drawn from $\sim\mathcal{N}(0,(9\times10^{-4})^2)$ 
\begin{equation}
    \mathcal{N}(\mu, \sigma^2) = \frac{1}{\sqrt{2\pi\sigma^2}}\exp{\left[-\frac{(x-\mu)^{2}}{2\sigma^{2}}\right]}
\end{equation}
where $\mu$ is the normal distribution's mean and $\sigma$ is its standard deviation. We then find a mass-optimal forced periodic trajectory for each particle (if a solution is feasible), treating the state of the particle as the initial condition to the trajectory. The fitness of each particle is then computed as
\begin{equation}
    F = \boldsymbol{\psi}^T \delta\mathbf{x}
    \label{eq:fitness_weight}
\end{equation}
where $\boldsymbol{\psi}$ is a vector of weights. If a trajectory is infeasible, $F$ is set to 0 for that particle. The states of the particles are then updated with 
\begin{equation}
    \delta\mathbf{x}_0^{k+1} = \beta (\mathbf{g} - \delta\mathbf{x}_0^k) + \alpha \mathbf{N}
    \label{eq:state_update}
\end{equation}
where $\mathbf{g}$ is the state of the most fit particle. $\beta \in (0,1)$ controls how a particle will be updated compared to the best known particle's state. The greater $\beta$ is, the faster the swarm will converge on the current best. $\alpha \in (0,1)$ controls how the particles explore the proximity of the search space. In our implementation, we set $\beta = 0.7$ and $\alpha = 0.5^k$ where $k$ is the index of the iteration. With these parameter selections, we balance exploration of the search space and fast convergence of the swarm. $\mathbf{N}$ is a vector of parameters chosen from a normal distribution at every iteration and is scaled to the size of the search space. Using our previous work on this topic\cite{Merrill2024} to inform the standard deviations of this vector, we find that
\begin{equation}
    \mathbf{N} =
    \begin{bmatrix}
      \mathcal{N}(0,(8\times10^{-3})^2) \\
      \mathcal{N}(0,(8\times10^{-4})^2)\\ 
      \mathcal{N}(0,(6\times10^{-3})^2)\\
      \mathcal{N}(0,(8\times10^{-3})^2)\\
      \mathcal{N}(0,(1\times10^{-2})^2)\\
      \mathcal{N}(0,(4\times10^{-3})^2)
    \end{bmatrix} 
\end{equation}
works well for this problem. This iteration process is repeated until our stopping criteria are reached. We choose to terminate the process when 20 iterations are completed and the maximum fitness has not increased or when $J_M/(u_{max}(t_f-t_0)) > 0.95$ (i.e., more than $95\%$ of the trajectory is spent thrusting). If the PSO has taken more than $\sim 80$ iterations but $J_M/(u_{max}(t_f-t_0)) < 0.85$, the fitness evaluation is switched to
\begin{equation}
    F = J_M
    \label{eq:fitness_dV}
\end{equation}
with the current set of particles so that particles that increase $J_M$ are explored more thoroughly. After a solution achieves $J_M/(u_{max}(t_f-t_0)) > 0.85$, the fitness evaluation is switched back to Equation \ref{eq:fitness_weight}. 

In this paper, the procedure described above is repeated for 
\begin{equation}
    \boldsymbol{\psi} = 
    \begin{bmatrix}
      \cos(\psi)\\
      \sin(\psi)\\
      \mathbf{0}_{4\times1}
    \end{bmatrix} 
\end{equation}
where $\psi$ is increased from $0$ to $2\pi$ in increments of $\pi/6$. A representative example of a PSO run's performance is shown in Figure \ref{fig:pso_example}. On the right plot in Figure \ref{fig:pso_example}, the darkness of the particles are scaled by their iteration number, where the solid particles are those from iteration 14 and the most transparent belong to the initial swarm. Across almost every iteration, only $\sim$ half of the particles converge to a solution. In this problem, a low convergence percentage indicates that the particles are either outside of the region of convergence entirely or that there is a dense swarm near the barrier between feasible and infeasible trajectories. The $0\%$ convergence shown at the first iteration of Figure \ref{fig:pso_example} is the former case, where the particles have been heavily dispersed, as $\alpha=1$ for that iteration and the randomness of particles is maximized. Beyond iteration 10, we observe that the average fitness of the converged particles approaches a steady state and rarely changes significantly.  


\section{Results}

\subsection{Mass-Optimal Forced Periodic Trajectories}

\begin{figure}[htbp]
    \centering
    \includegraphics[width=1.0\linewidth]{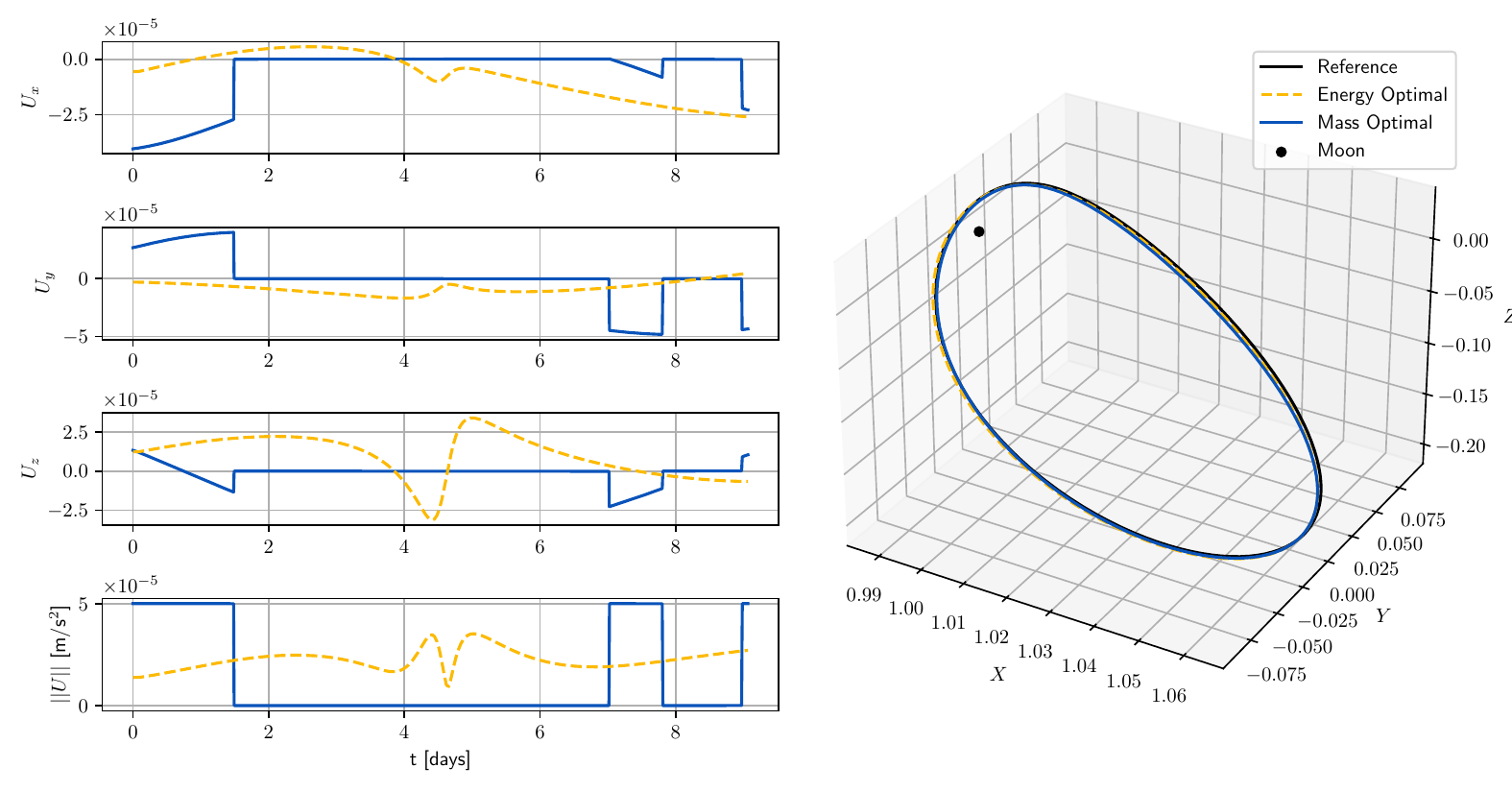}
    \caption{An example of energy-optimal and mass-optimal trajectories where the timing of the thrusting maxima are dissimilar between the two cases.}
    \label{fig:mass_dissim}
\end{figure}

\begin{figure}[htbp]
    \centering
    \includegraphics[width=1.0\linewidth]{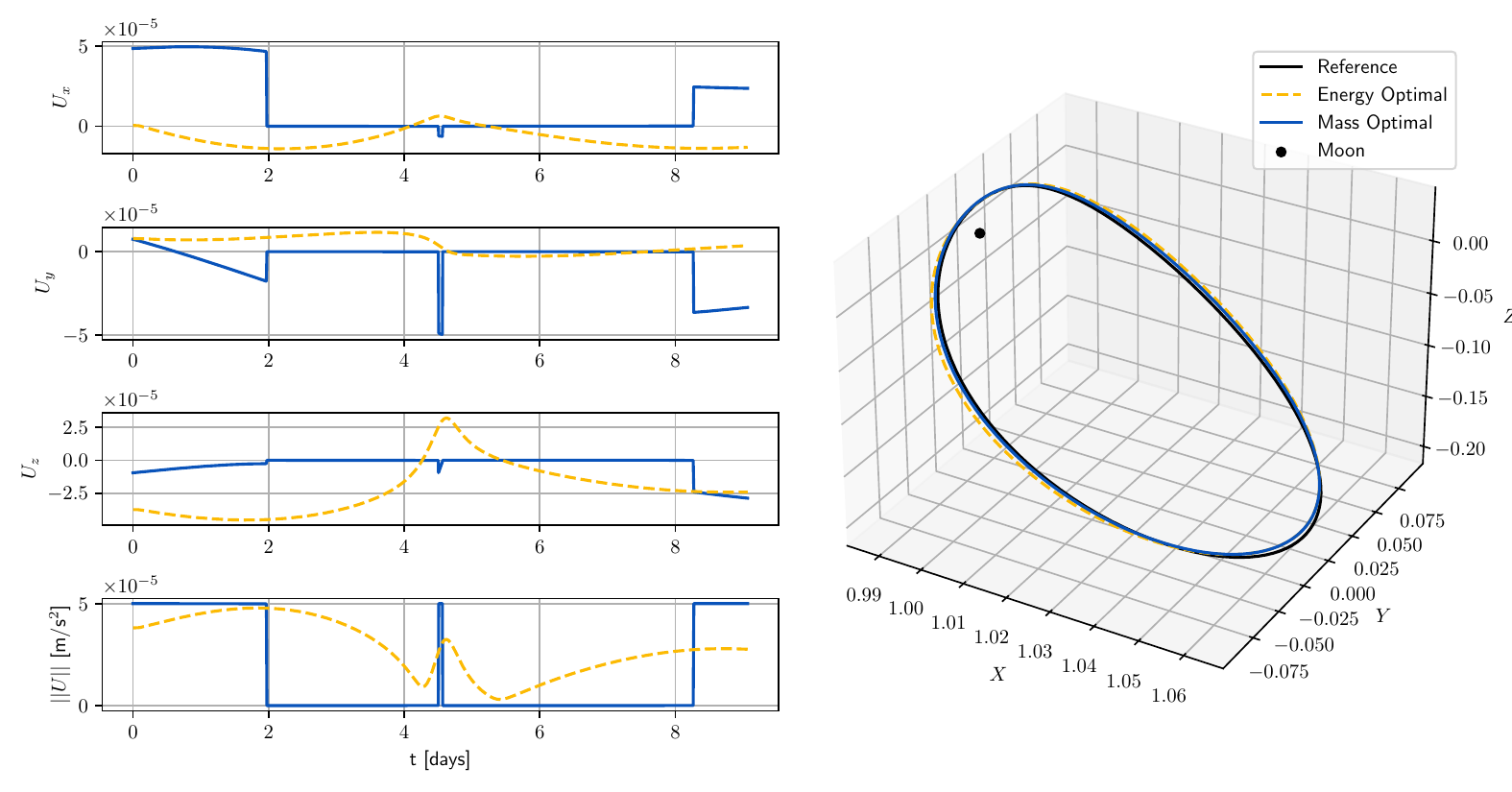}
    \caption{An example of energy-optimal and mass-optimal trajectories where the timing of the thrusting maxima are similar between the two cases.}
    \label{fig:mass_sim}
\end{figure}

Here, we provide two representative examples of forced periodic trajectories that exist in the proximity of our reference trajectory. Given the constraints in these problems, the maximum $\Delta V$ that may be expended in a single orbit is $\sim 39.3$ m/s. In the formulation of the mass-optimal problem, $\Delta V = J_M$. The mass-optimal trajectory in Figure \ref{fig:mass_dissim} has a $\Delta V = 10.2$ m/s per orbit and the mass-optimal trajectory in Figure \ref{fig:mass_sim} has a $\Delta V = 12.2$ m/s per orbit. 

In Figure \ref{fig:mass_dissim}, the mass-optimal trajectory features thrusting at the minima of the energy-optimal trajectory. This contrasts with Figure \ref{fig:mass_sim}, where the thrusting of the mass-optimal trajectory is well-aligned with the maxima of the energy-optimal trajectory. Among the sampled trajectories, there is not a clear predictor of the mass-optimal thrust timing that can be gained by the energy-optimal trajectory. In almost every case, the mass-optimal trajectory will feature an extended burn at the beginning of the trajectory and a short burn near the half-period of the trajectory (when the spacecraft reaches perilune). All mass-optimal trajectories feature characteristic bang-bang thrust profiles, where the thruster either fires at maximum thrust or not at all. The energy-optimal trajectories feature a smooth and continuous thrust profile which often changes most drastically in magnitude near perilune. Although a thrust limit is not applied to the energy-optimal trajectories, many of them feature maximum thrust magnitudes similar to that of the mass-optimal trajectories (as are shown in the examples here). 

\subsection{Reachable Set Comparison}

The linearized energy-limited reachable set is a series of six-dimensional ellipsoids centered around the states along the reference trajectory. When projected to a two-dimensional plane, this reachable set takes the form of a series of ellipses. We then interpolate between the edges of the ellipses to find a volume around the reference trajectory. The non-linear thrust-limited reachable set, by comparison, does not have a well-defined shape. To approximate this reachable set, we find mass-optimal forced periodic trajectories where the spacecraft thrusts for at least $95\%$ of the trajectory. 

\begin{figure}[ht]
    \centering
    \includegraphics[width=1.0\linewidth]{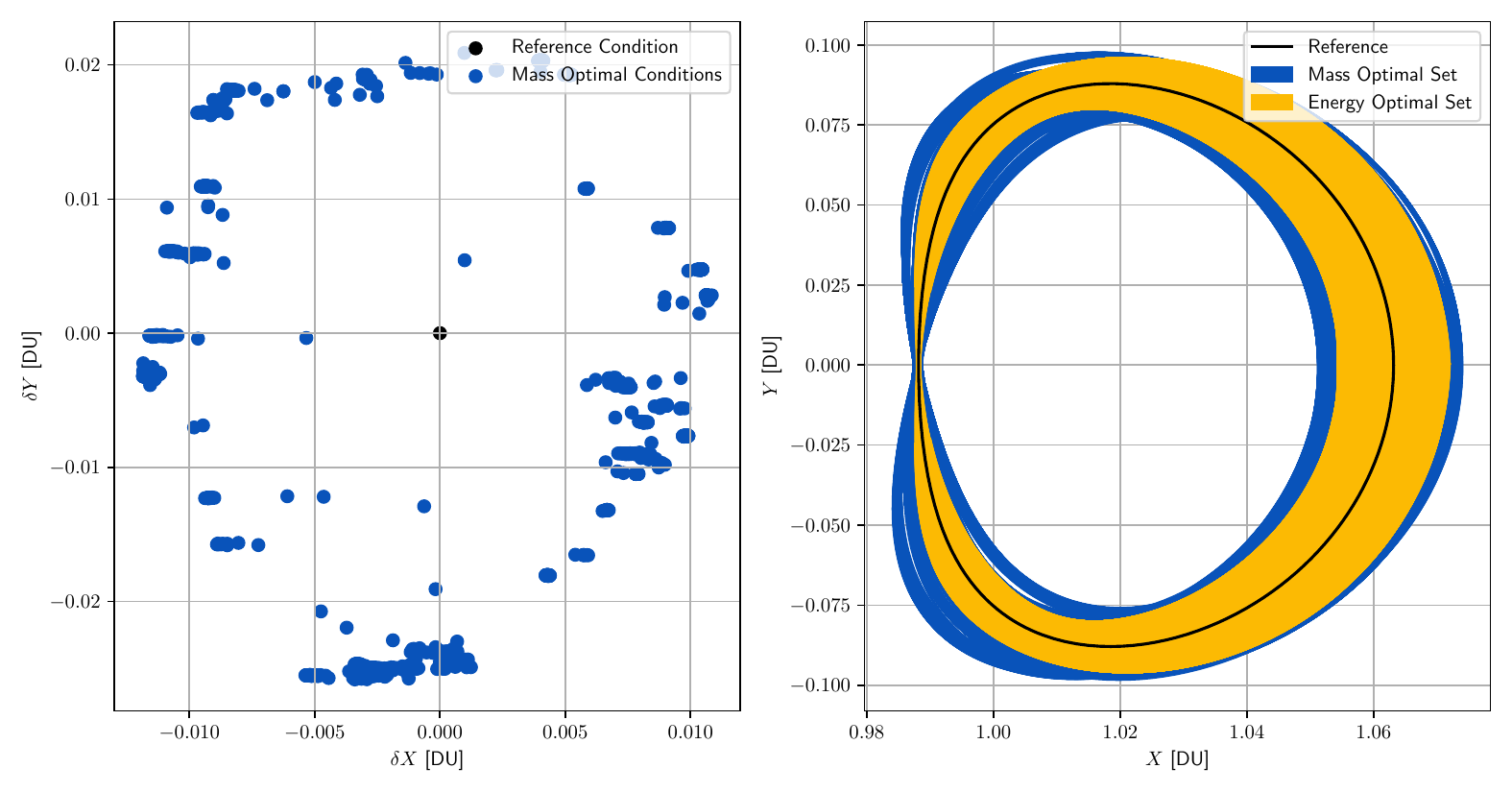}
    \caption{A comparison of the reachable sets for energy-optimal and mass-optimal trajectories in the $xy$-plane.}
    \label{fig:mass_reachability}
\end{figure}

Figure \ref{fig:mass_reachability} shows mass-optimal forced periodic trajectories that satisfy $J_M/(u_{max}(t_f-t_0)) > 0.95$. Although the full set of trajectories is incomplete due to the imperfect sampling, the energy-optimal reachable set based on linearization shares similarity to the full non-linear mass-optimal reachable set. It is important to note that if the same solution method is used to find an energy-optimal reachable set and then a mass-optimal reachable set, the energy-optimal reachable set would be a superset of the mass-optimal reachable set. However, because the two optimization methods that we use to obtain these sets use different assumptions and constraints, this is not the case. 


\section{Conclusions} \label{sec:conc}



The mass-optimal trajectories presented in this work are just the first of a variety of mass-optimal forced periodic trajectories. The generation code written on top of the ASSET tool is robust and capable of producing these trajectories reliably and quickly. Building continuation schemes for the computation of these trajectories and exploring their diversity is just one topic of our future work.

The generation of the thrust-limited reachable set via PSO is achieved for a subset of the reachable space in this paper. This is an expensive technique but does offer an improvement to random sampling of the space. The PSO approach is made more useful with the prior knowledge of the scale of the problem found by the previous linear analysis. Our PSO algorithm is also improved by initializing the particles around solutions to similar problems. 

We find that the reachable set of nonlinear mass-optimal thrust-limited trajectories is a superset of the linearized energy-optimal energy-limited trajectories. This is likely due to the alterations made to the dynamics via the linearization. We do note that these sets are very similar when projected in the $xy$-plane. The investigation of these sets projected into other spaces is the subject of future work.

\section*{Acknowledgements}
The authors thank Carrie Sandel for her support with setting up ASSET and her helpful comments. This material is based upon work supported by the Air Force Office of Scientific Research under award number FA9550-23-1-0665.



\bibliographystyle{AAS_publication}   
\bibliography{references}   


\end{document}